\documentclass[12pt]{amsart} 
\usepackage{amssymb,amsmath,graphicx,colordvi}

\newtheorem{theorem}{Theorem} 
\newtheorem{lemma}[theorem]{Lemma}
\newtheorem{corollary}[theorem]{Corollary}

\copyrightinfo{}{}

\title{Irrational proofs for three theorems of Stanley} 

\author{Matthias Beck}
\address{Department of Mathematics\\
         San Francisco State University\\
         San Francisco, CA 94132\\
         USA}
\email{beck@math.sfsu.edu}
\urladdr{http://math.sfsu.edu/beck}

\author{Frank Sottile}
\address{Department of Mathematics\\
         Texas A\&M University\\
         College Station, TX 77843-3368\\
         USA}
\email{sottile@math.tamu.edu}
\urladdr{http://www.math.tamu.edu/\~{}sottile}

\thanks{Work of Sottile is supported in part by NSF CAREER grant DMS-0134860}
\thanks{We thank an anonymous referee and Robin Chapman for pointing
  out flaws in an earlier version of this 
  article, Christian Haase for suggesting that we look for an irrational proof of
  Stanley's Monotonicity Theorem, and Richard Stanley for his comments on that proof.}

\renewcommand{\L}{{\mathcal L}}
\newcommand{\K}{{\mathcal K}}
\renewcommand{\P}{{\mathcal P}}
\newcommand{\Q}{{\mathcal Q}}

\newcommand{\R}{{\mathbb R}}
\newcommand{\N}{{\mathbb N}}
\newcommand{\Z}{{\mathbb Z}}
\newcommand{\m}{{\bf m}}
\newcommand{\s}{{\bf s}}
\newcommand{\sv}{{\bf s}}  
\renewcommand{\v}{{\bf v}}

\newcommand{\w}{{\bf w}}
\newcommand{\x}{{\bf x}}

\newcommand\cone{\operatorname{cone}} 
 
\begin{document}

\begin{abstract}
 We give new proofs of three theorems of Stanley on generating functions
 for the integer points in rational cones.
 The first 
 relates the rational generating function 
 $\sigma_{\v+\K} (\x)\ :=\ \sum_{ \m \in \left( \v+\K \right) \cap \Z^d } \x^\m$,
 where $\K$ is a rational cone and $\v \in \R^d$, with $\sigma_{ -\v+\K^\circ } (1/\x)$.
 The second theorem 
 asserts that the generating function $1 + \sum_{ n \ge 1 } L_\P (n) \, t^n$
 of the Ehrhart
 quasi-polynomial $L_\P (n) := \# \left( n\P \cap \Z^d \right)$ of a rational polytope
 $\P$ can be written as a rational function $\frac{ \nu_\P (t) }{ \left( 1-t \right)^{ \dim \P + 1} }$
 with \emph{nonnegative} numerator $\nu_\P$.
 The third theorem 
 asserts that if $\P\subseteq\Q$, then 
 $\nu_\P\leq\nu_\Q$.
 Our proofs are based on elementary 
 counting afforded by irrational
 decompositions of rational polyhedra.
\end{abstract}

\maketitle


\section{Introduction}

For us, a (\emph{convex}) \emph{rational polyhedron} $\P$ is the intersection of finitely many
half-spaces in $\R^d$, where each half-space has the form 
$\left\{ \left( x_1, x_2, \dots, x_d \right) \in \R^d \mid \ a_1 x_1 + a_2 x_2 + \dots + a_d x_d \le b \right\}$
for some integers $a_1, a_2, \dots, a_d, b$.
A \emph{rational cone} $\K$ is a rational polyhedron with a unique vertex at the origin. 
We are interested in the generating function
 \[
  \sigma_{\v+\K} (\x)\ :=\ \sum_{ \m \in \left( \v+\K \right) \cap \Z^d } \x^\m
 \]
for the integral points of the shifted (``affine") cone $\v+\K$ and its companion $\sigma_{\v+\K^\circ}(\x)$ for
the integral points of the (relative) interior $\K^\circ$ of $\K$. 
Here, $\x^\m$ denotes the product $x_1^{m_1}x_2^{m_2} \cdots x_d^{m_d}$. 
The function $\sigma_{\v+\K}$ (as well as $\sigma_{ \v+\K^\circ }$) is a rational function in
the variables $\x$. 
Stanley's \emph{Reciprocity Theorem}~\cite{stanleyreciprocity} relates the 
functions $\sigma_{\v+\K}$ and $\sigma_{ -\v+\K^\circ }$ for any $\v \in \R^d$. 
We abbreviate the vector $\left( 1/x_1, 1/x_2, \dots, 1/x_d \right)$ by 
$\frac{1}{\x}$.

\begin{theorem}[Stanley]\label{stanrec}
  Suppose that $\K$ is a rational cone and $\v \in \R^d$.
  Then, as rational functions, 
  $\sigma_{\v+\K} \left( \x \right) = (-1)^{ \dim \K } \, 
   \sigma_{-\v+\K^\circ}\left(\frac{1}{\x}\right)$.
\end{theorem}

There are proofs of Theorem~\ref{stanrec} which involve local cohomology in commutative
algebra~\cite[Section I.8]{stanleycombcommalg} and complex
analysis~\cite{stanleyreciprocity}. 
Many proofs, including ours, first prove it for the easy
case of simplicial cones, and then use a decomposition of $\K$ into simplicial cones
to deduce Theorem~\ref{stanrec}. 
This approach requires some device to handle the
subsequent overcounting of integral points that occurs as the cones in the decomposition
overlap along 
faces.
In other proofs, this device is either a shelling argument~\cite{ziegler}, or a
valuation (finitely additive measure)~\cite{klainrota}, or some other version of
inclusion-exclusion.
In contrast, our method of `irrational decomposition' requires no such device as the
proper faces of the cones we use contain no integral points.

We use the same construction to prove Stanley's \emph{Positivity Theorem}.
A  \emph{rational polytope} is a bounded rational polyhedron.
A rational polytope is \emph{integral} if its vertices lie in $\Z^d$.
For an integral polytope $\P \subset \R^d$, Ehrhart~\cite{ehrhartpolynomial}
showed that  the function 
\[
  L_\P (n)\ :=\ \# \left( n \P \cap \Z^d \right) 
\]
is a polynomial in the integer variable $n$. 
If the polytope $\P$ is only rational, then the function $L_\P(n)$ is a
\emph{quasi-polynomial}.
More precisely, let $p$ be a positive integer such that $p\P$ is integral.
Then there exist polynomials $f_0,f_1,\dotsc,f_{p-1}$ so that
\[
  L_\P (n)\ =\ f_{(n\, \mbox{\scriptsize mod}\, p)}(n)\,.
\]
(It is most efficient, but not necessary, to take the minimal such $p$.)

The generating function for $L_\P$ is a rational function with denominator
$(1-t^p)^{\dim\P+1}$ (see, for example, \cite[Chapter 4]{stanleyec1} or the proof we
give in Section~\ref{S:positivity}). 
But one can say  more~\cite{stanleydecomp}.

\begin{theorem}[Stanley]\label{stanpos}
 Suppose $\P$ is a rational $d$-polytope with $p\P$ integral and set
 \begin{equation}\label{Eq:Polytope-Hilbert}
    1 + \sum_{ n \ge 1 } L_\P (n) \, t^n\ =\ 
     \frac{ a_{(d+1)p-1} t^{(d+1)p-1} + a_{(d+1)p-2} t^{(d+1)p-2} + \dots + a_0 }{ (1-t^p)^{ d+1 } } \ .
 \end{equation}
 Then $a_0, a_1, \dots, a_{(d+1)p-1} \ge 0$.
\end{theorem}

Even more can be said.
Suppose that $\Q$ is a rational polytope containing $\P$ and that both $p\P$ and
$p\Q$ are integral. 
Supressing their dependence on $p$, let $\nu_\P$ and $\nu_\Q$ be the numerators of the
rational generating functions~\eqref{Eq:Polytope-Hilbert} for $\P$ and $\Q$,
respectively. 
We have $d=\dim\P<\dim\Q=e$ and so $\nu_\Q$ is the numerator of the rational
generating function for $L_\Q(n)$, which has denominator $(1-t^p)^e$, while
$\nu_\P$ is the numerator of the rational generating function for $L_\P(n)$,
which has denominator $(1-t^p)^e$.
Stanley's \emph{Monotonicity Theorem}~\cite{stanleymon} aserts that every coefficient of
$\nu_\Q$ dominates the corresponding coefficient of $\nu_\P$, that is,
$\nu_\P\leq\nu_\Q$. 

\begin{theorem}[Stanley]\label{stanineq}
 Suppose $\P\subseteq \Q$ are rational polytopes with $p\P$ and $p\Q$ integral.
 Then $\nu_\P \leq \nu_\Q$.
\end{theorem}

While Theorem \ref{stanrec} may seem unconnected to Theorems \ref{stanpos} and
\ref{stanineq},
they are related by a construction which---to the best of our
knowledge---is due to Ehrhart. 
Lift the  vertices $\v_1,\v_2,\dots,\v_m$ of a rational polytope $\P \subset
\R^d$ into $\R^{1+d}$, by adding 1
as their first coordinate, and let $p$ be a positive integer such that $p\P$ is
integral. 
Then 
\[
  \v_1'\ =\ \left( p, p\v_1 \right) ,\ \
  \v_2'\ =\ \left( p, p\v_2 \right) ,\ \
   \dots,\ \ \v_m' = \left(p , p\v_m \right)\ 
\]
are integral.
Now we define the \emph{cone over} $\P$ to be
\[
  \cone (\P)\ =\ \left\{ \lambda_1 \v_1' + \lambda_2 \v_2' + \dots + \lambda_m \v_m' 
   \mid \ \lambda_1, \lambda_2, \dots, \lambda_m \ge 0 \right\} \subset \R^{ 1+d } .
\]
We can recover our original polytope $\P$ (strictly speaking, the set 
 $ \left\{ (1,\x) \mid  \ \x \in \P \right\} $) by cutting $\cone (\P)$ with the
 hyperplane $x_0 = 1$.  
Cutting $\cone(\P)$  with the hyperplane $x_0 = 2$, we obtain a copy of $2 \P$, cutting 
with $x_0 = 3$ gives a copy of $3 \P$, etc. 
Hence
\[
  \sigma_{ \cone (\P) } \left( x_0, x_1, \dots, x_d \right) 
  = 1 + \sum_{ n \ge 1 } \sigma_{ n \P } \left(x_1,\ldots,x_d\right) x_0^n \ . 
\]
Since 
$\sigma_{n\P}\left(1,1,\dots,1\right) = \#\left( n \P \cap \Z^d \right)$, we obtain
 \[ 
    \sigma_{ \cone (\P) } \left( t, 1, 1, \dots, 1 \right) 
     \ =\ 1 + \sum_{ n \ge 1 } L_\P (n) \, t^n \ .
\]

A nice application of Theorem \ref{stanrec} 
is the following reciprocity theorem, which was conjectured (and partially
proved) by Ehrhart \cite{ehrhart2} and proved 
by Macdonald \cite{macdonald}. 

\begin{corollary}[Ehrhart-Macdonald]\label{emrec}
The quasi-polynomials $L_\P$ and $L_{\P^\circ}$ satisfy
\[
  L_\P (-t) = (-1)^{\dim \P} L_{\P^\circ} (t) \ .
\]
\end{corollary} 

As with Theorem \ref{stanrec}, most proofs of Theorem~\ref{stanpos} use shellings of a
polyhedron or finite additive measures (see, e.g.,
\cite{ehrhartbook,macdonald,mcmullenreciprocity}). The only exceptions we are aware of
are proofs via complex analysis (see, e.g., \cite{stanleyreciprocity}) and commutative
algebra (see, e.g., \cite[Section I.8]{stanleycombcommalg}). 
We feel that no existing proof is as elementary as the one we give.

We remark that the same technique gives a similarly elementary and
subtraction-free proof of Brion's Theorem~\cite{Br88}.
This proof will appear in~\cite{BHS}.


\section{Stanley's Reciprocity Theorem for cones} 

Any cone has a triangulation into \emph{simplicial cones} which are cones with a
minimal number of boundary hyperplanes (see, e.g., \cite{leetriangulations}).  
This is the starting point for our proof, which differs from other proofs that
use such a  decomposition.
The decomposition that we use is, from the view of integer points,
non-overlapping, and thus we need only apply elementary (as in
elementary--school) counting arguments, sidestepping any hint of
inclusion-exclusion. 

\begin{proof}[Irrational Proof of Theorem~$\ref{stanrec}$]
 Triangulate $\K$ into simplicial rational cones $\K_1,\K_2,\dots,\K_n$, all having the
 same dimension as $\K$. 
 Now there exists a vector $\s \in \R^d$ such that 
\begin{equation}\label{shiftedconeinterior}
  \left( \v+\K^\circ \right) \cap \Z^d\ =\ (\s+\K) \cap \Z^d
\end{equation}
and
\begin{equation}\label{triangconesboundary}
  \partial \left( \pm\s + \K_j \right) \cap \Z^d\ =\ 
   \emptyset \qquad \text{ for all } j=1, \dots, m\, .
\end{equation}
In fact, $\s$ may be any vector in the relative interior of some cone
$\v+\K_i$ for which $\s-\v$ is short enough such that 
\eqref{shiftedconeinterior} holds. 

This means, in particular, that there are no lattice points on the
boundary of $\v + \K$, and so \eqref{shiftedconeinterior} implies 
$ (-\v+\K) \cap \Z^d = (-\s+\K) \cap \Z^d $.
Furthermore, because of \eqref{triangconesboundary},
\[
  \sigma_{-\v+\K} \left( \x \right)
  \ =\ \sigma_{-\s+\K} \left( \x \right)
  \ =\ \sum_{ j=1 }^{ m } \sigma_{-\s+\K_j} \left( \x \right)
\]
and
\[
  \sigma_{ \v+\K^\circ } \left( \x \right)
  \ =\  \sigma_{ \s+\K } \left( \x \right)
  \ =\ \sum_{ j=1 }^{ m } \sigma_{ \s+\K_j } \left( \x \right)\,.
\]
The result now follows from reciprocity for simplicial cones, which is Lemma
\ref{simpleconereciprocity} below.
\end{proof}

Despite our title, the vector $\s-\v$ need not be irrational  as any
short rational vector will do.

\begin{lemma}\label{simpleconereciprocity}
Fix linearly independent vectors $\w_1, \w_2, \dots, \w_d \in \Z^d$, and let
\[
  \K\ =\ \left\{ \lambda_1 \w_1 + \lambda_2 \w_2 + \dots + \lambda_d \w_d \mid
     \ \lambda_1, \dots, \lambda_d \ge 0 \right\}\, .
\]
Then for those $\sv \in \R^d$ for which the boundary of the shifted simplicial cone 
$\sv + \K$ contains no integer point, 
\[
  \sigma_{ \sv + \K } \left({\textstyle \frac 1 \x }\right)\ =\
     (-1)^d \, \sigma_{ -\sv + \K } \left( \x \right)\, .
\]
\end{lemma}

As in Theorem \ref{stanrec}, the reciprocity identity is one of rational functions. 
In the course of the proof, we will show that  $\sigma_{ \sv + \K }$ is indeed a
rational function for $\sv \in \R^d$.

\begin{proof}
If we tile the cone $\sv+\K$ with $\N\{\w_1,\w_2,\dotsc,\w_d\}$--translates of the
half-open parallelepiped $\sv+\P$, where 
\[
  \P\ :=\ \left\{ \lambda_1 \w_1 + \lambda_2 \w_2 + \dots + \lambda_d \w_d 
     \mid \ 0 \le \lambda_1, \lambda_2, \dots, \lambda_d < 1 \right\}\, ,
\]
then we can express $\sigma_{\sv+\K}$ using geometric series
 \begin{equation}\label{simpleconegenfct}
    \sigma_{\sv+\K} (\x)\ =\ \frac{\sigma_{\sv+\P}(\x)}{\left(1-\x^{\w_1}\right) 
       \left( 1 - \x^{ \w_2 } \right) \cdots \left( 1 - \x^{ \w_d } \right) } \ .
 \end{equation}
(This proves that $\sigma_{ \sv + \K }$ is a rational function.)
Similarly, 
\[
  \sigma_{-\sv+\K} (\x)\ =\ \frac{ \sigma_{-\sv+\P} (\x) }{ \left(1 -\x^{ \w_1 }\right) 
    \left( 1 - \x^{ \w_2 } \right) \cdots \left( 1 - \x^{ \w_d } \right) } \ ,
\]
so we only need to relate the parallelepipeds $\sv+\P$ and $-\sv+\P$. 
By assumption, $\sv+\P$ contains no integer points on its boundary,
and so we may replace $\P$ by its closure.
Note that $\P=\w_1 + \w_2 + \dots + \w_d -\P$, so we have 
the identity
 \begin{equation}\label{parallelid}
    \sv+\P\ =\ -(-\sv+\P) + \w_1 + \w_2 + \dots + \w_d \ .
\end{equation}

In terms of generating functions, \eqref{parallelid} implies that
\[
  \sigma_{\sv+\P} \left( \x \right)\ =\
  \sigma_{ -\sv+\P } \left({\textstyle \frac 1 \x} \right) 
   \x^{ \w_1 } \x^{ \w_2 } \cdots \x^{ \w_d } \ ,
\]
whence
\begin{eqnarray*}
  \sigma_{\sv+\K} \left({\textstyle \frac 1 \x} \right) &=& \frac{ \sigma_{ \sv+\P } \left( \frac 1 \x \right) }{ \left( 1 - \x^{ -\w_1 } \right) \left( 1 - \x^{ -\w_2 } \right) \cdots \left( 1 - \x^{ -\w_d } \right) } \\
                 &=& \frac{ \sigma_{ -\sv+\P } (\x) \, \x^{ -\w_1 } \x^{ -\w_2 } \cdots \x^{ -\w_d } }{ \left( 1 - \x^{ -\w_1 } \right) \left( 1 - \x^{ -\w_2 } \right) \cdots \left( 1 - \x^{ -\w_d } \right) } \\
                 &=& \frac{ \sigma_{ -\sv+\P } (\x) }{ \left( \x^{ \w_1 } - 1 \right) \left( \x^{ \w_2 } - 1 \right) \cdots \left( \x^{ \w_d } - 1 \right) } \\
                 &=& (-1)^d \frac{ \sigma_{ -\sv+\P } (\x) }{ \left( 1 - \x^{ \w_1 } \right) \left( 1 - \x^{ \w_2 } \right) \cdots \left( 1 - \x^{ \w_d } \right) } \\
                 &=& (-1)^d \, \sigma_{ -\sv+\K } (\x) .
\end{eqnarray*}
\end{proof}

Lemma~\ref{simpleconereciprocity} is essentially due to Ehrhart.
The new idea here is our `irrational' decomposition.


\section{Stanley's Positivity and Monotonicity Theorems for Ehrhart
  polynomials}\label{S:positivity} 

\begin{proof}[Irrational Proof of Theorem~$\ref{stanpos}$]
As before, triangulate $\cone(\P) \subset \R^{ d+1 } $ into simple rational cones 
 $\K_1, \K_2, \dots, \K_m$, each of whose generators are among the generators
 $(p,p\v_i)$ of $\cone(\P)$.
(Such a triangulation always exists; see, e.g., \cite{leetriangulations}.)
Again there exists a vector $\sv \in \R^{d+1}$ such that
\[
  \cone(\P) \cap \Z^d\ =\ (\sv+\cone(\P)) \cap \Z^d
\]
and no facet of any cone $\sv+\K_i$ contains any integral points. 
Thus every integral point in $\sv+\cone(\P)$ belongs to exactly one simplicial cone
$\sv+\K_j$, and we have 
\[
  \cone(\P) \cap \Z^d\ =\
  (\sv+\cone(\P)) \cap \Z^d\ =\ 
  \bigcup_{ j=1 }^m \left( \left( \sv+\K_j \right) \cap \Z^d \right)\, ,
\]
and this union is \emph{disjoint}. 
We obtain the identity of generating functions, 
\[
  \sigma_{ \cone(\P) } \left( \x \right)\ =\ 
  \sum_{ j=1 }^{ m } \sigma_{ \sv+\K_j } \left( \x \right) .
\]
But now we recall from the introduction that
\[
  1 + \sum_{ n \ge 1 } L_\P (n) \, t^n\ =\
    \sigma_{ \cone(\P) } \left( t, 1, 1, \dots, 1 \right)\  =\  
  \sum_{ j=1 }^{ m } \sigma_{ \sv+\K_j } \left( t, 1, 1, \dots, 1 \right)\, .
\]
 So it suffices to show that the rational generating functions 
 $\sigma_{ \sv+\K_j } \left( t, 1, 1, \dots, 1 \right)$ for the \emph{simplicial}
 cones $\sv+\K_j$ have nonnegative numerators and denominators of the form
 $(1-t^p)^{d+1}$. 

 In this case, the cone $\sv+\K_j$ has integral generators of the form
 $\w_i=(p,p\v_i)$, for some vertices $\v_1,\dotsc,\v_{d+1}$ of the polytope $\P$,
 where $p$ is a positive integer such that $p\P$ is integral.
 Substituting $( t, 1, 1, \dots, 1)$ into the concrete form of the rational
 generating function \eqref{simpleconegenfct}, gives denominator 
 $(1-t^p)^{d+1}$ and numerator the generating function for the integer points
 in the parallelepiped which is generated by $\w_1,\dotsc,\w_{d+1}$ and shifted
 by $\sv$, where the coefficient $a_i$ of $t^i$ counts points with first
 coordinate $i$. 
\end{proof}

\begin{proof}[Irrational Proof of Theorem~$\ref{stanineq}$]
 Suppose first that $\dim\P=\dim\Q$.
 As in the previous proof, suppose that  $\K_1, \K_2, \dots, \K_m$ triangulate
 $\cone(\P)$ into simplicial rational cones, each of whose generators are among the
 generators $(p,p\v_i)$ of $\cone(\P)$ \cite{leetriangulations}.
 We may extend this to a triangulation  $\K_1, \K_2, \dots, \K_l$ of 
 $\cone(\Q)$, where the additional simplicial cones have generators from the given
 generators  $(p,p\v_i)$ of $\cone(\P)$ and $(p,p\w_i)$ of $\cone(\Q)$.
 The generators of each cone $\K_i$ and the irrational shift vector $\sv$
 together give a parallelepiped with no lattice points on its boundary,
 and the coefficient of $t^j$ in $\nu_\P$ is the number of integer points with last
 coordinate $j$ in the union of these parallelepipeds for $\K_1,\dotsc,\K_m$. 
 The result follows as the coefficient of $t^j$ in $\nu_\Q$ is the number of
 integer points with last coordinate $j$ in the parallelepipeds for
 $\K_1,\dotsc,\K_l$, and $m<l$.  

 If however, $\dim\P<\dim_Q$, then the triangulation  $\K_1, \K_2, \dots, \K_m$ 
 of $\cone(\P)$ extends to a triangulation $\L_1,\L_2,\dots,\L_l$ of
 $\cone(\Q)$, where now the simplicial cones $\K_i$ are $d$-faces of the simplicial cones
 $\L_j$.
 Note that the irrational decomposition $\sv+\L_j$, $j=1,\dotsc,l$ restricts to
 an irrational decomposition of $\cone(\P)$ given by some vector
 $\sv'\in\R\cdot \cone(\P)$.
 Moreover, for every $i=1,\dotsc,m$ there is a unique $a(i)$ with 
 $1\leq a(i)\leq l$ such that $\sv'+\K_i\subset\sv+\L_{a(i)}$. 
 The same is true for the parallelepipeds generated by the vectors $(p,\v)$ along
 the rays of these cones, and also for their shifts by $\sv'$ and $\sv$.
 Then the result follows by the same argument as before once we interpret the
 coefficients of $t^j$ in $\nu_\P$ and $\nu_\Q$ as the number of points with
 second coordinate $j$ in the union of these parallelepipeds.
\end{proof}

\def\cprime{$'$} \def\cprime{$'$}
\providecommand{\bysame}{\leavevmode\hbox to3em{\hrulefill}\thinspace}
\providecommand{\MR}{\relax\ifhmode\unskip\space\fi MR }
\providecommand{\MRhref}[2]{%
  \href{http://www.ams.org/mathscinet-getitem?mr=#1}{#2}
}
\providecommand{\href}[2]{#2}

\end{document}